\documentclass{icmart}

\contact[thomas.schick@mathematik.uni-goettingen.de]{Mathematisches Institut,
  Georg-August-Universit\"at G\"ottingen, Bunsenstr. 3, 37073 G\"ottingen,
  Germany}

\usepackage{color}
\usepackage{epsfig}
\usepackage{amscd}

\theoremstyle{plain}
\newtheorem{theorem}{Theorem}[section]

\newtheorem{proposition}[theorem]{Proposition}
\newtheorem{conjecture}[theorem]{Conjecture}

\theoremstyle{definition}
\newtheorem{definition}[theorem]{Definition}
\newtheorem{example}[theorem]{Example}

\newtheorem{convention}[theorem]{Convention}

\newtheorem{remark}[theorem]{Remark}
\newtheorem{question}[theorem]{Question}

\newcommand{\reals}{\mathbb{R}}
\newcommand{\complexs}{\mathbb{C}}
\newcommand{\naturals}{\mathbb{N}}
\newcommand{\integers}{\mathbb{Z}}

\DeclareMathOperator{\id}{id}
\newcommand{\boundary}[1]{\partial#1}
\newcommand{\boundedops}{\mathcal{B}}
\newcommand{\compacts}{\mathcal{K}}

\newcommand{\norm}[1]{\left\lVert#1\right\rVert}

\newcommand{\tensor}{\otimes}
\newcommand{\into}{\hookrightarrow}

\newcommand{\iso}{\cong}
\newcommand{\disjointunion}{\amalg}

\DeclareMathOperator{\Riemplus}{Riem^+}
\DeclareMathOperator{\eucl}{eucl}
\DeclareMathOperator{\lf}{lf}
\DeclareMathOperator{\Pos}{Pos}
\DeclareMathOperator{\Diff}{Diffeo}

\DeclareMathOperator{\supp}{supp}   
\DeclareMathOperator{\im}{im}      
\DeclareMathOperator{\vol}{vol}    

\DeclareMathOperator{\spin}{spin}

\DeclareMathOperator{\ind}{ind}
\DeclareMathOperator{\sgn}{sgn}
\DeclareMathOperator{\scal}{scal}  

\newcommand{\forget}[1]{}

\newcommand{\innerprod}[1]{\langle #1 \rangle}

{\catcode`@=11\global\let\c@equation=\c@theorem}
\renewcommand{\theequation}{\thetheorem}

\allowdisplaybreaks[2]

\title{The topology of positive scalar curvature}

\author[Thomas Schick]{Thomas Schick\thanks{supported by the Courant
    Research Center "`Higher order structures in Mathematics"' of     Georg-August-Universit\"at G\"ottingen}}

\begin{document}

\begin{abstract}
  Given a smooth closed manifold $M$ we study the space of
  Riemannian  metrics
  of positive scalar curvature on $M$. A long-standing question is: when is
  this space non-empty (i.e.~when does $M$ admit a metric of positive
  scalar curvature)? More generally: what is the topology of this space?
  For example, what are its homotopy groups?

  Higher index theory of the Dirac operator is the basic tool to address these
  questions. This has seen tremendous development in recent years, and in
this survey we will discuss some of the most pertinent examples. 

In particular, we will show how advancements of \emph{large scale index
  theory} (also called \emph{coarse index theory}) give
rise to new 
types of obstructions, and provide the tools
for a systematic study of the existence and classification problem via the
K-theory of $C^*$-algebras. This is part of a program  ``mapping the
topology of positive scalar curvature to analysis''.

In addition, we will show how advanced surgery theory and smoothing theory
can be used to construct the first elements of infinite order in
the $k$-th homotopy groups of the space of metrics of positive scalar
curvature for arbitrarily large $k$.
Moreover, these examples 
are the first ones which remain non-trivial in the moduli space of such
metrics. 
\end{abstract}



\maketitle

\section{Introduction}

  One of the fundamental questions at the interface of geometry and topology
  concerns the relation between local geometry and global topology. 

More specifically, given a compact smooth manifold $M$ without boundary, what
are the possibilities for Riemannian metrics on $M$? Even more specifically,
can we find a metric of positive scalar curvature on $M$ and if yes, what does
the space of such metrics look like?

Recall the following definition of the scalar curvature
function.
\begin{definition}
  Given an $n$-dimensional smooth Riemannian manifold $(M,g)$, the \emph{scalar
    curvature} at $x$ describes the volume expansion of small balls around $x$
  via 
  \begin{equation*}
    \frac{\vol(B_\epsilon(M,x))}{\vol(B_\epsilon(\reals^n,0))} = 1-
    \frac{\scal(x)}{6(n+2)} \epsilon^2 + O(\epsilon^4),
  \end{equation*}
  compare \cite[0.60]{Besse}. In particular, this means that if $\scal(x)>0$
  then geodesic balls around $x$ for small radius have smaller volume than the
  comparison balls in Euclidean space. 
  Of course, alternatively, $\scal(x)$  can be defined as a second
  contraction of the Riemannian curvature operator at $x$.
\end{definition}

The most important tool to investigate these questions goes back to Erwin
Schr\"odinger \cite{Schroedinger}, rediscovered by Andr\'e Lichnerowicz
\cite{Lichnerowicz}: If $M$ has positive 
scalar curvature and a spin structure then the Dirac operator on $M$ is
invertible. This forces its index (which is the super-dimension of the null
space) to vanish.

Recall that a spin structure is a (global) differential geometric datum for a
Riemannian manifold $M$ which allows us to construct a specific Riemannian
vector bundle $S$, the spinor bundle, together with a specific differential
operator of order $1$, the Dirac operator $D$ (compare
e.g.~\cite{LawsonMichelsohn} for a nice introduction). 

On the other hand
non-vanishing of the index follows from index theorems, giving rise to
powerful \emph{obstructions} to positive scalar
curvature. For example, the Atiyah-Singer index theorem says that
$\ind(D)=\hat A(M)$, 
where the $\hat A$-genus is a fundamental differential topological invariant (not
depending on the metric!). 

The most intriguing  question around this method to rule out positive scalar
curvature asks to what extent a sophisticated refinement of $\ind(D)$, the
\emph{Rosenberg index} $\alpha^{\reals}(M)$ which takes values in the
K-theory of the (real)
$C^*$-algebra of the fundamental group $\Gamma$ of $M$, is the only
obstruction. This is the content of the (stable) Gromov-Lawson-Rosenberg
conjecture.

\begin{conjecture}\label{conj:GLR}
  Let $M$ be a connected closed spin manifold of dimension $\ge 5$. The
  \emph{Gromov-Lawson-Rosenberg conjecture} asserts that $M$ admits a metric
  with positive scalar curvature if and only if $\alpha^{\reals}(M)=0 \in
  KO_*(C_{\reals}^*\pi_1(M))$. 

  The \emph{stable Gromov-Lawson-Rosenberg conjecture} claims that
  $\alpha^{\reals}(M)=0$ if and only if there is $k\in\mathbb{N}$ such that
  $M\times 
  B^k$ admits a metric with positive scalar curvature. Here, $B$ is any
  so-called \emph{Bott manifold} i.e.~a
  simply connected $8$-dimensional spin manifold with $\hat A(M)=1$.
\end{conjecture}

To put the stable version in context: given two closed manifolds $M,B$ such
that $M$ admits a metric of positive scalar curvature, so does $M\times B$,
simply using the product of a sufficiently scaled metric on $M$ with any
metric on $B$. Therefore the unstable Gromov-Lawson-Rosenberg conjecture
implies the stable one.

Recall here that for a discrete group $\Gamma$ the \emph{maximal group
  $C^*$-algebra} 
  $C^*_{\max}\Gamma$ is defined as the completion of the group ring
  $\complexs[\Gamma]$ with respect to the maximal possible $C^*$-norm on
  $\complexs[\Gamma]$, and the
  \emph{reduced group $C^*$-algebra} 
  $C^*_{\text{r}}\Gamma$ is defined as the norm closure of $\complexs[\Gamma]$,
  embedded in $\boundedops(l^2(\Gamma))$ via the regular representation. The
  \emph{real} group
  $C^*$-algebras $C^*_{{\reals},r}\Gamma$ and $C^*_{\reals,max}\Gamma$ replace
  $\complexs$ by $\reals$ throughout. Using them gives more information,
  necessary in the 
  Gromov-Lawson-Rosenberg conjecture. In this survey, for simplicity we will
  not discuss them but concentrate on the complex versions. For the
  Rosenberg index one can use them all, where a priori
  $\alpha_{\max}(M)\in K_*(C^*_{\max}\pi_1(M))$ is stronger than
  $\alpha_{\text{r}}(M)\in K_*(C^*_{\text{r}}\pi_1(M))$. The \emph{Baum-Connes isomorphism
    conjecture}, compare \cite{Baum-Connes-Higson}, predicts the calculation
  of   $K_*(C^*_{\text{r}}\Gamma)$ in terms of the
  equivariant K-homology of a suitable classifying space. The \emph{strong
    Novikov conjecture} predicts that this equivariant K-homology at least
  embeds.

In celebrated work Stephan Stolz \cite{Stolz1,StolzTrieste} has
established the following two partial positive results.

\begin{theorem}\label{theo:Stolz}
  The Gromov-Lawson-Rosenberg conjecture is true for manifolds with trivial
  fundamental group. In other words, if $M$ is a closed connected spin
  manifold of 
  dimension $\ge 5$ with trivial fundamental group, then $M$ admits a
  Riemannian metric with positive scalar 
  curvature if and only if $\alpha^{\reals}(M)=0$.

  More generally, if $\pi_1(M)$ satisfies the strong Novikov conjecture then
  the stable Gromov-Lawson-Rosenberg conjecture is true for $M$.
\end{theorem}

On the other hand, recall the counterexamples of
\cite{Schick_counterex, DwyerSchickStolz} which show that the
unstable Gromov-Lawson-Rosenberg conjecture is not always true.
\begin{theorem}\label{theo:counterex}
  For $5\le n\le 8$ there exist closed spin manifolds $M^n$ of dimension $n$
  such that   $\alpha(M^n)=0$, but such that $M^n$ does not admit a metric with
  positive scalar curvature.
\end{theorem}

The manifolds $M_n$ can be constructed with fundamental groups
$\mathbb{Z}^{n-1} \times \mathbb{Z}/3\mathbb{Z}$ or with appropriately chosen
torsion-free fundamental group, but for example not with a
free abelian fundamental group \cite{JoachimSchick}. It remains one of the
most intriguing open 
questions whether the Gromov-Lawson-Rosenberg conjecture is true for
all $n$-dimensional manifolds with fundamental group
$(\mathbb{Z}/3\mathbb{Z})^n$.

The obstructions used in the counterexamples of Theorem \ref{theo:counterex}
are not based on index theory, but on the minimal 
hypersurface method of Richard Schoen and Shing-Tung Yau \cite{SchoenYau}
which we will not discuss further in this survey. 

As a companion to the Gromov-Lawson-Rosenberg conjecture we suggest a
slightly weaker conjecture about
the strength of the Rosenberg index:
\begin{conjecture}\label{conj:Rosenberg_power}
  Let $M$ be a closed spin manifold. Every obstruction to positive scalar
  curvature for manifolds of dimension $\ge 5$ which is based on
  index theory of Dirac operators can be read off the Rosenberg index
  $\alpha^{\reals}(M)\in KO_*(C_{\reals}^*\pi_1(M))$.
\end{conjecture}

This is vague because the statement ``based on index theory of Dirac operators''
certainly leaves room for interpretation. 

By Stolz' Theorem \ref{theo:Stolz}, Conjecture \ref{conj:Rosenberg_power}
follows from the strong Novikov conjecture. On
 the other hand,
every index theoretic obstruction which is not (yet) understood in terms of
the Rosenberg index is particularly interesting. After all, it is a potential
starting point to obtain counterexamples to the strong Novikov conjecture.

Around this question we discuss the following results
\cite{HankeSchick1,HankeSchick2,HankeKotschickRoeSchick}.

\begin{theorem}\label{theo:enlarg}
  Let $M$ be an area-enlargeable spin manifold (which implies by the work of
  Mikhail Gromov and Blaine Lawson \cite{GromovLawson} that 
  $M$ does not admit a metric of positive scalar curvature).

  Then $\alpha_{\max}(M) \ne 0 \in K_*(C^*_{\max}\pi_1(M))$.

  If $M$ is even (length)-enlargeable, then $\alpha_{\text{r}}(M)\ne 0 \in
  K_*(C^*_{\text{r}}\pi_1(M))$. 
\end{theorem}

Recall that a closed $n$-dimensional manifold $M$ is called \emph{enlargeable}
if it admits a 
sequence of coverings $M_i\to M$ which come with compactly
supported maps $f_i\colon M_i\to S^n$ of non-zero degree but such that
$\sup_{x\in M_i} \norm{D_xf_i}$ tends to $0$ as $i\to\infty$. It is
\emph{area-enlargeable} if the same holds with $\norm{\Lambda^2D_xf_i}$ (a
weaker condition). For
the definition of the norms, we use a fixed metric on $M$ and its pull-backs
to $M_i$ and a fixed metric on $S^n$.

As a potential counterexample to Conjecture \ref{conj:Rosenberg_power} we
describe a \emph{codimension-$2$} obstruction to
positive scalar curvature (in a special form introduced by Mikhail Gromov and
Blaine Lawson in \cite[Theorem
7.5]{GromovLawson}) which is based on index theory of the Dirac
operator, but which so far is not known to be encompassed by the
Rosenberg index.

\begin{theorem} (compare \cite[Section 4]{HankePapeSchick}).\label{theo:codim2}
    Let $M$ be a closed connected spin manifold with vanishing second homotopy
    group. Assume that $N\subset M$ is a smooth submanifold
    of codimension $2$ with trivial normal bundle and such that the inclusion
    induces an injection on the level of fundamental groups $\pi_1(N)\into
    \pi_1(M)$. Finally, assume that the Rosenberg index of the Dirac operator
    on the 
    submanifold $N$
     does not vanish: $0\ne \alpha(N)\in K_*(C^*\pi_1N)$.

 Then $M$ does not admit a Riemannian metric with positive scalar curvature.
\end{theorem}

(Secondary)
index invariants of the Dirac operator 
can be used in the \emph{classification} of metrics of positive scalar
curvature, if applied to appropriately constructed
examples. ``Classification'' means in particular to understand how many
deformation classes of metrics of positive scalar curvature a given manifold
carries, or more generally what the topology of the space of such metrics looks
like. 

A promising tool to systematically study the existence and classification
problem is the ``Stolz positive scalar curvature long exact
sequence''. It has the form
\begin{equation*}
\cdots  \to R_{n+1}(\pi_1(M)) \to \Pos_n(M) \to \Omega^{\spin}_n(M)\to
R_n(\pi_1(M))\to \cdots 
\end{equation*}
Here the group we would like to understand is $\Pos_n(M)$, the structure
group of metrics of positive scalar curvature 
(on the spin manifold $M$ and related manifolds, and modulo a suitable bordism
relation). The group $\Omega^{\spin}_n(M)$ 
is the usual spin bordism group from algebraic topology, which is very well
understood. Finally, $R_*(\pi_1(M))$ indeed depends only on the fundamental
group of the manifold in question.
Note that this positive scalar curvature sequence is very similar in spirit to
the surgery exact sequence coming up in the classification of manifolds.

Unfortunately we have not yet been able to fully compute all the
terms in this exact sequence, even for the simplest possible case of trivial
fundamental group. However, a lot of
information can be gained using index theory by mapping out to more manageable
targets. Here, we refer in particular to \cite{piazzaschickrho},
joint with Paolo Piazza, where we construct a commutative diagram of
maps, using \emph{large scale index theory}, to 
the K-theory sequence of associated $C^*$-algebras
\begin{equation*}
\cdots\to  K_{n+1}(C_{\text{r}}^*\Gamma) \to K_{n+1}(D^*\tilde
M^\Gamma)\to K_n(M)\xrightarrow{\alpha} K_n(C_{\text{r}}^*\Gamma)\to \cdots
\end{equation*}
We again abbreviate $\Gamma=\pi_1(M)$. This sequence was introduced by Nigel Higson
and John Roe \cite{HR1} and called there the \emph{analytic surgery exact
  sequence}. A lot is known about this K-theory
sequence: $K_n(M)$ is just the usual topological K-homology of $M$ (an
important generalized homology theory). Moreover, $\alpha$ is the Baum-Connes 
assembly map.

A good deal of this survey will discuss the \emph{large scale index theory}
underlying the constructions. In particular, we will explain two
primary and secondary index theorems which play key roles in the application
of this theory: 
\begin{itemize}
\item A vanishing theorem for the large scale index of the Dirac operator under
  partial positivity of the scalar curvature, Theorem
  \ref{theo:partialvanish}. 
\item A higher secondary index theorem, that shows how the large
  scale index of a manifold with boundary (which has positive scalar curvature
  near 
  the boundary) determines a structure invariant of the boundary's
  metric of positive scalar curvature.
\end{itemize}

A fundamental problem in the use of (higher) index theory centers around the
question whether there is a difference between \emph{topological information}
(which typically can be computed much more systematically) and analytical
information (which has the desired geometric consequences but often is hard to
compute). This is answered (conjecturally) by the strong Novikov
conjecture. This explains why these
conjectures play such a central role in higher index theory. The appropriate
version for large scale index theory is the \emph{coarse Baum-Connes
  conjecture} (here, ``with coefficients'', for details compare Section
\ref{sec:coarseBC}).
 
\begin{conjecture}
    Given a locally compact metric space $X$ of bounded geometry and an auxiliary coefficient
    $C^*$-algebra $A$, then in the composition
  \begin{equation*}
    K_*^{\lf}(X;A) \to KX_*(X;A) \to K_*(C^*(X;A))
  \end{equation*}
  the second map is an isomorphism. Here, $K_*^{\lf}(X)$ is the topologists'
  locally finite K-homology (a generalized homology theory) of the space $X$,
  and $K^{\lf}_*(X;A)$ is a version with coefficients, still a generalized
  cohomology theory. Finally $KX_*(X;A)$, the coarse K-homology, is a variant
  which depends only on the large scale geometry of $X$. 
\end{conjecture}

This conjecture has many concrete applications. It implies
 the strong Novikov conjecture. However, I expect that counterexamples to
 these conjectures eventually
will be found.

Concerning the classification question mentioned above, in the last part of
the survey we will discuss a new construction method, based on
advanced surgery theory and smoothing theory in the topology of manifolds. 

\begin{definition}
 We define $\Riemplus(M)$ to be the space of Riemannian metrics of positive
 scalar curvature on $M$, an infinite dimensional manifold.
\end{definition}

The
main result is that $\pi_k(\Riemplus(M))$ is very often non-trivial, even its
image in the moduli space of such metrics remains non-trivial.

More precisely, we have the following theorem, derived in joint work with
Bernhard Hanke and Wolfgang Steimle (compare \cite[Theorem 1.1]{HSS}).
\begin{theorem}\label{theo:large_pi}
  For every $k\in\naturals$, there is $n_k\in \naturals$ such that, whenever
  $M$ is a connected closed spin manifold with a metric $g_0$ of positive scalar
  curvature and with $\dim(M)>n_k$ and $k+\dim(M)+1\equiv 0 \pmod 8$,
  then $\pi_k(\Riemplus(M),g_0)$ contains an element of infinite order. 

  If $M$ is a homology sphere, then the image of this element in
  $\pi_k(\Riemplus(M)/\Diff_{x_0}(M))$ also has infinite order, where the
  diffeomorphism group acts by pullback.
\end{theorem}

The second part of the theorem implies that the examples constructed do not
rely on the homotopy properties of the diffeomorphism group of $M$. This is in
contrast to all previous 
known cases, compare in particular \cite{Hitchin,CrowleySchick}.

Here, $\Diff_{x_0}(M)$ is the subgroup of the full diffeomorphism group
consisting of diffeomorphisms of $M$ which fix the point $x_0\in M$ and whose
differential at $x_0$ is the identity. It is much more reasonable to use this
subgroup instead of the full diffeomorphism group, because it ensures that the
moduli 
space $\Riemplus(M)/\Diff_{x_0}(M)$ remains an infinite dimensional
manifold, instead of producing a very singular space.


\begin{remark}
  Most of the results mentioned so far display also how poorly the
  topology of positive scalar curvature is understood: the method relies on
  the index theory of the Dirac operator and the Schr\"odinger-Lichnerowicz
  formula. This is quite a miraculous relation which certainly is very
  helpful. But it requires the presence of a spin structure. Manifolds without
  spin structure (and where not even the universal covering admits a spin
  structure) a priori shouldn't be very different from manifolds with spin
  structure, i.e.~one would expect that many of them do not admit a metric of
  positive scalar curvature. But until now we have almost no tools to decide
  this (apart from the minimal surface method, which is only established in small
  dimensions). 
 
  Almost any progress in this direction would be a real
  breakthrough.
\end{remark}

\section{Index theory and obstructions to positive scalar curvature}
\label{sec:index_i_principle}
The underlying principle how scalar curvature is coupled to the Dirac operator
comes from a formula of Schr\"odinger \cite{Schroedinger}, rediscovered and
first applied by 
Lichnerowicz \cite{Lichnerowicz}. The starting point is a spin manifold
$(M,g)$, with spinor bundle $S$ and
\emph{Dirac operator} $D$. Schr\"odinger's formula says
\begin{equation*}
  D^2 = \nabla^*\nabla + \frac{\scal}{4}.
\end{equation*}
The first term on the right is the ``rough Laplacian'', by definition a
non-negative
unbounded operator on the $L^2$-sections of $S$. The second term stands for
point-wise multiplication with the scalar curvature function. If
the scalar curvature is bounded below by a positive number (called ``uniformly
positive'' later), this is a positive operator and
hence $D$ is invertible.

Let us recall the basics of the index theory of the Dirac operator, formulated
in the language of operator algebras and K-theory. This is the most
convenient setup for the generalizations we have in mind.

We start with a very brief introduction to the K-theory of $C^*$-algebras.
\begin{enumerate}
\item The assignment $A\mapsto K_*(A)$ is a functor from the category of
  $C^*$-algebras to   the category of graded abelian groups. 
\item We can (for a unital $C^*$-algebra $A$) define $K_0(A)$ as the group of
  equivalence classes of projectors in $A$ and the matrix algebras $M_n(A)$.
\item We can (for a unital $C^*$-algebra $A$) define $K_1(A)$ as the group of
  equivalence classes of invertible elements in $A$ and $M_n(A)$.
\item There is a natural Bott periodicity isomorphism $K_n(A)\to
  K_{n+2}(A)$. 
\item For each short exact sequence of $C^*$-algebras $0\to I\to A\to Q\to 0$
  there is naturally associated a long exact sequence in K-theory
  \begin{equation*}
    \cdots \to K_{n+1}(Q)\xrightarrow{\delta} K_n(I)\to K_n(A)\to
    K_n(Q)\to\cdots .
  \end{equation*}
\item One can generalize K-theory for real and graded $C^*$-algebras. In the
  former case Bott periodicity has period $8$.
\item On can use extra symmetries based on Clifford algebras to give
  descriptions of $K_n(A)$ which are adapted to the treatment of
  $n$-dimensional spin manifolds.
\end{enumerate}

On an even dimensional manifold, the spinor bundle canonically
splits
into $S=S^+\oplus S^-$ and the Dirac operator is odd, i.e.~has the form $D=
{\begin{pmatrix}
  0 & D^-\\ D^+ & 0
\end{pmatrix}}$. Let $\chi\colon\reals\to\reals$ be any continuous function
which is odd (i.e.~$\chi(-x)=\chi(x)$ for $x\in\reals$) and with
$\lim_{x\to\infty}\chi(x)=1$. Functional calculus allows to define $\chi(D)$,
which is an odd bounded operator acting on $L^2(S)$. Choosing an
isometry $U\colon L^2(S^-)\to L^2(S^+)$ we form
$U\chi(D)^+ \in \boundedops:=\boundedops(L^2(S^+))$, the $C^*$-algebra of all
bounded operators on $L^2(S^+)$. If $M$ is \emph{compact}, ellipticity of $D$
implies that $U\chi(D)^+$ is invertible modulo the ideal
$\compacts:=\compacts(L^2(S^+))$ of compact operators. The short exact
sequence of 
$C^*$-algebras $0\to\compacts\to \boundedops\to
\boundedops/\compacts\to 0$
gives rise to a long exact K-theory
sequence and the relevant piece of this
exact sequence for us is
\begin{equation}\label{eq:basicKseq}
\cdots \to K_1(\boundedops)  \to K_1(\boundedops/\compacts) \xrightarrow{\delta}
  K_0(\compacts) \to\cdots
\end{equation}
As invertible elements in $A$ represent classes in $K_1(A)$, the above
spectral considerations yield a class $[U\chi(D)^+] \in
K_1(\boundedops/\compacts)$ 
and we define the index to be $\ind(D):=\delta(U\chi(D)^+) \in
K_0(\compacts)$.  

Of course, for the compact operators, $K_0(\compacts)$ is isomorphic
to $\integers$, generated by any rank $1$ projector in
$\compacts$. In our case
$\delta(U\chi(D)^+)$ is represented by the projector onto the kernel of $D^+$
minus the projector onto the kernel of $D^-$, so that we
arrive at the usual $\ind(D):=\dim\ker(D^+)-\dim\ker(D^-)$.

Analysing the situation more closely, the additional geometric information of
positive scalar curvature implies invertibility of $D$ which translates to the
fact that $U\chi(D)^+$ is invertible already in $\boundedops$ and
then $[U\chi(D)^+]\in
K_1(\boundedops/\compacts)$ in the sequence
\eqref{eq:basicKseq} has a lift to
$K_1(\boundedops)$. Exactness implies $\ind(D)=0$.

A second main ingredient concerning the index of the Dirac operator is the
Atiyah-Singer {index theorem} \cite{ASIII}. A
priori $\ind(D)$ (like the operator $D$) depends on the Riemannian metric on
$M$. However, the {index theorem} expresses it in terms which are
independent of the metric. More specifically, $\ind(D)= \hat A(M)$, the $\hat
A$-genus of $M$.

This result has vast generalizations in many directions. A very important
one (introduced by Jonathan Rosenberg \cite{Rosenberg_PSC_NC_I}) modifies the
Dirac operator by ``twisting'' with a smooth flat bundle $E$ of 
(finitely generated projective) modules over an auxiliary $C^*$-algebra
$A$. One then obtains an index in 
$K_*(A)$. Indeed, the construction is pretty much the same as above, with the
important innovation that one replaces the scalars $\complexs$ by the more
interesting $C^*$-algebra $A$, as detailed in Section
\ref{sec:subs-theory-manif}.



The second generalization works for a non-compact manifold $X$. In this case the classical Fredholm
property of the Dirac operator fails. To overcome this,  
\emph{large scale index theory}, synonymously called \emph{coarse index
  theory} is developing. Again this is based on $C^*$-techniques and
pioneered by John Roe
\cite{Roe_index_coarse}. It is
tailor-made for the non-compact setting. One obtains an index in the K-theory
of the Roe algebra $C^*(X;A)$. In $C^*$-algebras, positivity implies
invertibility, which finally implies that all the generalized indices
vanish if one starts with a metric of uniformly positive scalar curvature.

The general pattern (from the point of view adopted in this article) of index
theory is the following:
\begin{enumerate}
\item The geometry of the manifold $M$ produces an interesting operator $D$.
\item This operators defines an element in an operator algebra $A$, which
  depends on the precise context.
\item The operator satisfies a Fredholm condition, which means it is
  invertible module an ideal $I$ of the algebra $A$, 
  again depending on the context.
\item The algebras in question are $C^*$-algebras.  This implies that
  ``positivity'' of elements is defined, and moreover ``positivity'' implies
  invertibility.
\item A very special additional geometric input implies positivity and hence
  honest invertibility of our operator. For us, this special context will be
  the fact that we deal with a metric of uniformly positive scalar curvature.
\item Indeed, any element which is invertible in $A$ modulo an ideal $I$
  defines an element in $K_{n+1}(A/I)$, where $n=\dim(M)$.  (Instead of getting
  $K_1$, the fact that we deal with the Dirac
  operator of an $n$-dimensional manifold produces additional symmetries
  (related to actions of the Clifford algebra $Cl_n$) which give rise to the
  element in $K_{n+1}(A/I)$.)
\item 
 We interpret the class defined by the Dirac operator
  as a fundamental class $[M]\in K_{n+1}(A/I)$. Homotopy invariance of
  K-theory implies that $[M]$ does not depend on the full geometric data
  which goes in the construction of the operator $D$, but only on the
  topology of $M$.
\item The K-theory exact sequence of the extension $0\to
  I\to A\to A/I\to 0$ 
  contains the boundary map $\delta$. We call the image of
  $[M]$ under $\delta$ the index
  \begin{equation*}
    \delta\colon K_{n+1}(A/I)\to K_n(I); \quad [M]\mapsto \ind(D).
  \end{equation*}
  Note that the degree arises from additional dimension-dependent symmetries
  which we do not discuss in this survey.
\item The additional \emph{geometric} positivity assumption (uniformly
  positive scalar curvature) which 
  implies invertibility already in $A$, gives rise to a canonical lift of
  $[M]$ to an element $\rho(M,g)\in K_{n+1}(A)$. Because of this, we think of
  $K_{n+1}(A)$ as a \emph{structure group} and $\rho(M,g)$ is a \emph{structure
    class}. It   contains information about the underlying geometry.
\end{enumerate}

Indeed, we want to advocate here the idea that the setup just described has
quite a number of different manifestations, depending on the situation at
hand. It can be adapted in rather flexible ways. The next section treats one example.


\section{Large scale index theory}
\label{sec:subs-theory-manif}

We describe ``large scale index theory'' for a complete Riemannian manifold of
positive dimension. 

Therefore, let $(M,g)$ be such a complete Riemannian manifold. Fix a Hermitian
vector bundle $E\to M$ of positive dimension. We first describe the operator
algebras which are relevant. They are all defined as norm-closed subalgebras
of $\boundedops(L^2(M;E))$. 

\begin{definition}\label{def:Roe_alg}
  We need the following concepts.
  \begin{itemize}
  \item An operator $T\colon L^2(M;E)\to L^2(M;E)$ has \emph{finite
      propagation} (namely $\le R$) if $\phi T\psi=0$ whenever $\phi,\psi\in
    C_c(M)$
    are compactly supported continuous functions whose supports have distance
    at least $R$.

    Here, we think of $\phi$ also as bounded operator on $L^2(M;E)$, acting by
    point-wise multiplication.
  \item $T$ as above is called \emph{locally compact} if $\phi T$ and $T\phi$
    are compact operators whenever $\phi\in C_c(M)$.
  \item $T$ is called \emph{pseudolocal} if $\phi T \psi$ is compact whenever
    $\phi,\psi\in C_c(M)$ with disjoint supports, i.e.~such that $\phi\psi=0$.
  \item The \emph{Roe algebra} $C^*(M)$ is defined as the norm closure of the
    algebra of all bounded finite propagation operators which are locally
    compact. It is an ideal in the \emph{structure algebra} $D^*(M)$ which is
    defined as the closure of the algebra of finite propagation pseudolocal
    operators.
  \item Assume that a discrete group $\Gamma$ acts by isometries on
    $M$. Requiring  in the above definitions that the finite propagation
    operators are in   addition $\Gamma$-equivariant and then completing, we
    obtain the pair 
    $C^*(M)^\Gamma \subset D^*(M)^\Gamma$.  
  \end{itemize}
\end{definition}

\begin{remark}\label{rem:stabil}
  For technical reasons, one actually should replace the bundle $E$
  by the bundle $E\tensor l^2(\naturals)$  whose fibers are separable Hilbert
  spaces (or in the equivariant case by $E\tensor l^2(\naturals)\tensor
  l^2(\Gamma)$). Via the embedding $\boundedops(L^2(M;E))\to 
  \boundedops(L^2(M;E\tensor l^2(\naturals)))$ implicitly we think of
  operators on $L^2(M;E)$ as operators in the bigger algebra without
  mentioning this. Using the larger bundle guarantees functoriality
  and independence on $E$, implicit in our notion $D^*(M)$.
\end{remark}

  Let $A$ be an auxiliary $C^*$-algebra. Typical examples arise from a
  discrete group $\Gamma$, namely $C^*_{\max}\Gamma$ or $C^*_{\text{r}}\Gamma$.
  An important role is then played by smooth bundles $E$ over $M$  with fibers
  finitely 
  generated projective $A$-modules. These inherit fiberwise $A$-valued inner
  products (or more precisely Hilbert $A$-module structures in the sense of
  \cite{Lance}). Integrating the fiberwise inner product then also defines a
  Hilbert $A$-module structure on the space of continuous compactly supported
  sections of $E$. By completion, we obtain the Hilbert $A$-module
  $L^2(M;E)$. The bounded, 
  adjointable, $A$-linear operators on $L^2(M;E)$ form the Banach algebra
  $\boundedops(L^2(M;A))$. The ideal of \emph{$A$-compact operators} is defined
  as the norm closure of the ideal generated by operators of the form
  $s\mapsto v\cdot \innerprod{s,w}; L^2(M;E)\to L^2(M;E)$, with $v,w\in
  L^2(M;E)$.

  We now define $C^*(M;A)$ and $D^*(M;A)$ mimicking Definition
  \ref{def:Roe_alg}, but replacing the Hilbert space concepts by the Hilbert
  $A$-module concepts throughout. In particular, we use $A$-compact operators
  instead of compact operators. Also, the stabilization as discussed in Remark
  \ref{rem:stabil} is replaced suitably. 


\begin{example}\label{def:MishchenkoLine}
  Given a connected manifold $M$ with fundamental group $\Gamma$, there is a
  canonical $C^*\Gamma$-module bundle, the \emph{Mishchenko bundle}
  $\mathcal{L}$. It is the flat bundle associated to the left
  multiplication 
  action of $\Gamma$ on $C^*\Gamma$, where we treat $C^*\Gamma$ as the free
  right $C^*\Gamma$-module of rank $1$, i.e.~$\mathcal{L}=\tilde
  M\times_\Gamma C^*\Gamma$.
\end{example}

The construction of the large scale index is now  based on
 two principles.
\begin{proposition}\label{prop:D_prop}
  Let $M$ be a complete Riemannian spin manifold, $A$
  an auxiliary $C^*$-algebra and
  $E\to M$ a smooth bundle of (finitely generated projective) $A$-modules with
  compatible connection. We form the twisted Dirac operator $D_E$ as
  an unbounded operator on the Hilbert $A$-module of $L^2$-spinors on $M$ with
  values in $E$.

  For the operator $D_E$, there exists a functional calculus. In particular,
  we can form $f(D_E)$ for $f\colon \reals\to \reals$ a continuous function
  which vanishes at $\infty$ or which has limits as $t\to\pm\infty$. Moreover,
  $f(D_E)$ depends only on the restriction of $f$ to the spectrum of $D_E$.
  Then
  \begin{enumerate}
  \item if $f$ has a compactly supported (distributional) Fourier transform
    then $f(D_E)$ has finite propagation.
  \item if $f$ vanishes at infinity, then $f(D_E)$ is locally compact; if
    $f(t)$ converges as $t\to \pm\infty$ then $f(D_E)$ is at least
    pseudolocal. 
  \end{enumerate}
  The first property is a rather direct consequence of unit propagation speed
  for the fundamental solution of the heat equation. The second one is an
  incarnation of ellipticity and local elliptic regularity.   
\end{proposition}

These results are well known and have been used a lot in the literature
(compare in particular \cite{Roe_index_coarse}),
indeed they form the basis of ``large scale index theory''. For the very
general case needed in the proposition (with coefficients, arbitrary
complete $M$), a complete proof is given in \cite{HankePapeSchick}.

The construction of the large scale index is now rather straight-forward:
\begin{enumerate}
\item Take any continuous function $\psi\colon \reals\to \reals$ (later
  assumed to be odd) 
 with $\lim_{t\to\infty}\psi(t) =1$. Then Proposition
  \ref{prop:D_prop} implies that $\psi(D_E)$ belongs to $D^*(M;A)$. Even
  better, $1-\psi^2$ vanishes at infinity, so that $1-\psi(D_E)^2$ belongs to
  $C^*(M;A)$.
\item By the principles listed at the end of Section
  \ref{sec:index_i_principle}, $\psi(D_E)$ gives rise to an element $[M;E]$ in
  $K_{n+1}(D^*(M;A)/C^*(M;A))$. (Here, we again avoid the discussion of the
  additional symmetries which raise the index by $n$). Homotopy invariance of
  $C^*$-algebra K-theory implies that this element is independent of the
  choice of $\psi$ and depends only on the large scale features of the metric
  on $M$. 
\item We define the large scale index (or synonymously ``coarse index'')
  $\ind(D_E)\in K_n(C^*(M;A))$ as the image of $[M;E]$ under the boundary map
  of the long exact K-theory sequence.
\item If $M$ has uniformly positive scalar curvature and $E$ is flat, the
  Lichnerowicz-Weitzenb\"ock formula implies that $0$ is not in the spectrum
  of $D_E$. Then we can choose a function $\psi$ which is equal to $-1$ on the
  negative part of the spectrum of $D_E$ and equal to $+1$ on the positive
  part of the spectrum of $D_E$, so that $1-\psi^2$ vanishes on the spectrum
  of $D_E$, i.e.~$\psi^2(D_E)=1$.
  
  This means that $[M;E]$ lifts in a canonical way to $\rho(D;E)\in
  K_{n+1}(D^*(M;A))$ (this class depends on the metric of positive scalar
  curvature) and it implies that $\ind(D_E)=0$.

  
\item A special feature is that $K_{n+1}(D^*(M;A)/C^*(M;A))$ indeed is
  homological in nature: it is canonically isomorphic to the locally finite
  K-homology $K^{\lf}_n(M;A)$, satisfying the Eilenberg-Steenrod axioms of a
  (locally finite) generalized homology theory.
\end{enumerate}

\begin{example}
  If we apply this construction to a closed $n$-dimensional spin manifold
  $M$ and the Mishchenko bundle $\mathcal{L}$ on $M$, we obtain
$    \ind(D_{\mathcal{L}}) \in K_{n}(C^*(M;C^*\Gamma))$.

  However, there is a canonical isomorphism
  $K_*(C^*(M;C^*\Gamma)) \iso K_*(C^*\Gamma)$. 
  Using this isomorphism, the Rosenberg index mentioned
  above is exactly $\ind(D_{\mathcal{L}})$:
  \begin{equation*}
    \alpha(M) =\ind(D_{\mathcal{L}}) \in K_n(C^*\Gamma) \iso K_n(C^*(M;C^*\Gamma)).
  \end{equation*}
\end{example}

The reduced $C^*$-algebra $C^*_{\text{r}}\Gamma$ of a discrete group is a canonical
construction which captures many features of the group $\Gamma$,
e.g.~concerning its representation theory. However, it is very rigid. In
particular, it is not functorial: a homomorphism $\Gamma_1\to \Gamma_2$ will
in general not induce a homomorphism $C^*_{\text{r}}\Gamma_1\to
C^*_{\text{r}}\Gamma_2$. As a consequence it is very hard to find homomorphisms
out of $C^*_{\text{r}}\Gamma$ and also out of $K_*(C^*_{\text{r}}\Gamma)$. 

Coarse geometry, however, immediately provides such a homomorphism
(which allows one to detect elements in $K_*(C^*_{\text{r}}\Gamma)$).
This is based on a simple calculation:
 If a discrete group $\Gamma$ isometrically acts freely and cocompactly
  on a metric space $X$, then $C^*X^\Gamma$ is isomorphic to
  $C^*_{\text{r}}(\Gamma)\tensor \compacts$. ``Forgetting equivariance'' therefore
  gives the composite homomorphism 
  \begin{equation*}
  C^*_{\text{r}}\Gamma\into  C^*_{\text{r}}\Gamma\tensor\compacts \iso C^*X^\Gamma\into C^*X
  .
  \end{equation*}
  The induced map in K-theory allows one to detect elements in
  $K_*(C^*_{\text{r}}\Gamma)$ using large scale index theory, as we will show in one
  case in Section \ref{sec:enlarg}.

\section{The coarse Baum-Connes conjecture}
\label{sec:coarseBC}
Being the home of important index invariants, it is very important to be able
to compute the K-theory of the Roe algebras $C^*(M;A)$ for arbitrary complete
manifolds $M$ and coefficient $C^*$-algebras $A$. It turns out that
there are quite a number of tools to do this. Even better, at least
conjecturally there is a purely homological answer to this task.

Let us start with the three most important computational tools.
\begin{enumerate}
\item $K_*(C^*(M;A))$ is invariant under \emph{coarse homotopy}, compare
  \cite{HigsonRoe_homotopy}.
\item There are powerful \emph{vanishing theorems} for $K_*(C^*(M;A))$. An
  important one is valid if $M$ is \emph{flasque} \cite[Proposition
  9.4]{Roe_index_coarse}. This means that 
  $M$ admits a shift map $f\colon M\to M$ such that, on the one hand, $f$ is
  uniformly close to the identity (i.e.~there is a constant $C$ such that
  $d(f(x),x)<C$ for all $x\in M$). On the other hand, $f$ moves everything to
  infinity in the sense that for each compact subset $K$ of $M$, $\im(f^k)\cap
  K = \emptyset$ for all sufficiently large iterations $f^k$ of $f$.

\item The group $K_*(C^*(M;A))$ satisfies a Mayer-Vietoris principle. For
  this, one 
  needs a  \emph{coarsely excisive} decomposition $M=M_1\cup M_2$, which means
  that the 
  intersection $M_0:=M_1\cap M_2$ captures all the large scale features of the
  relation between $M_1$ and $M_2$. The technical definition is that for
  each 
  $R>0$ there is an $S>0$ such that the $S$-neighborhood of $M_1\cap M_2$
  contains the intersection of the $R$-neighborhoods of $M_1$ and $M_2$.

  In this situation, there is a long exact Mayer-Vietoris sequence (compare
  \cite{HigsonRoeYu,Siegel}) 
  \begin{multline*}
  \cdots\to    K_i(C^*(M_1;A))\oplus K_i(C^*(M_2;A))\to K_i(C^*(M;A)) \to\\
    K_{i-1}(C^*(M_0;A)) \to K_i(C^*(M_1;A))\oplus K_i(C^*(M_2;A))\to\cdots
  \end{multline*}
\end{enumerate}

One of the powerful principles for the K-theory of $C^*$-algebras is
their close relation to purely topological quantities via isomorphism
conjectures. Most prominent here is the Baum-Connes conjecture for
the computation of $K_*(C_{\text{r}}^*\Gamma)$. The properties of $K_*(C^*(M;A))$ listed
above indicate that a similar ``topological expression'' should be
possible here. Indeed, we have the \emph{coarse Baum-Connes conjecture (with
  coefficients)} \cite[Conjecture 8.2]{Roe_index_coarse}, verified in many
cases.

\begin{conjecture}
  Given a metric space $X$ of bounded geometry, in the composition
  \begin{equation*}
    K_*^{\lf}(X;A) \to KX_*(X;A) \to K_*(C^*(X;A))
  \end{equation*}
  the second map is an isomorphism.
\end{conjecture}

Here $K_*^{\lf}(X)$ is the usual locally finite K-homology of the space $X$,
defined analytically as $K_{*+1}(D^*X/C^*X)$, 
  and as we saw above it is no problem to introduce as coefficients a $C^*$-algebra $A$. The \emph{coarse K-homology $KX_*$} is
  obtained as the limit of $K_*^{\lf}(|\mathfrak{U}_i|)$, where the
  $\mathfrak{U}_i$ form a sequence of coverings of $X$ which become
  coarser as $i\to \infty$. Here $|\mathfrak{U}_i|$ is the geometric
  realization of the associated \v{C}ech simplicial complex. If $X$ is uniformly locally contractible, e.g.~if
  $X$ is the universal covering of a closed non-positively curved manifold,
  then
  the ``coarsening map'' $K_*^{\lf}(X;A)\to KX_*(X;A)$ is an isomorphism.

Recall that (in the context of large
scale geometry) ``bounded geometry'' means that $X$ contains a discrete
subset $T$ 
such that
on the one hand $T$ coarsely fills the space (i.e.~there is an $R>0$ such that the
  $R$-neighborhood of $T$ is all of $X$), but on the other hand $T$ is
  uniformly discrete (i.e.~for each $R>0$ the number of elements of
  $T$ contained in any $R$-ball is uniformly bounded from above).

The coarse Baum-Connes conjecture has a number of important consequences. Most
notably, there is a principle of descent \cite[Section 5]{Roe_index_coarse}
that uses the close relation 
between $C^*_{\text{r}}\Gamma$ and $C^*X$ for any metric space $X$ on which
$\Gamma$ acts properly and cocompactly. The principle of descent asserts that
if such a space 
satisfies the coarse Baum-Connes conjecture, then the strong Novikov
conjecture for $\Gamma$ is true.


On the other hand, the ``bounded geometry'' condition of the coarse
Baum-Connes conjecture is indispensable. Guoliang Yu has constructed a metric
space which is a disjoint union of (scaled) spheres of growing dimension for
which the analysis of the Dirac operator shows easily that the coarse assembly
map is not injective. Despite its simplicity, this example remains
intriguing. It is important to understand this better and to construct other
examples. We believe that the coarse
Baum-Connes conjecture will not be satisfied in full generality.

\section{Enlargeability and index}
\label{sec:enlarg}
Let $M$ be an (area)-enlargeable closed spin manifold. Recall that this
means that $M$ comes with a sequence of (not necessarily compact) coverings
$M_i$ with compactly supported maps of non-zero 
degree $f_i\colon M_i\to S^n$ which are arbitrarily (area) contracting.

Mikhail Gromov and Blaine Lawson show in \cite{GromovLawson} that an
enlargeable spin manifold does not admit a
metric of positive scalar curvature. 
Theorem \ref{theo:enlarg} verifies Conjecture \ref{conj:Rosenberg_power} for
this ``enlargeability obstruction'', i.e.~shows that the
Rosenberg index is non-zero in this situation. This was achieved  in
\cite{HankeSchick2} by refining the construction of Gromov and Lawson  as
follows: 
\begin{enumerate}
\item One constructs vector bundles $E_i$ on $M_i$ of small
  curvature which represent interesting K-theory classes.
\item Next one produces
  associated bundles $M(E_i)$ on $M$ (using a kind of push-down).
 If $M_i\to M$ is a
  finite covering, this is finite dimensional. If the
  covering $M_i\to M$ is infinite, 
  we canonically obtain an associated ``structure $C^*$-algebra''
  $C_i$ such that $M(E_i)$ is a Hilbert $C_i$-module bundle with finitely
  generated projective fiber.
\item The crucial step is the construction of a bundle $E:=\prod_{i}
  M(E_i)/\bigoplus_i M(E_i)$ which becomes a \emph{flat} Hilbert $A$-module
  bundle where $A=\prod_i C_i/\bigoplus_i C_i$. Being flat, this corresponds
  to a representation of $\pi_1(M)$.
\item As $E$ is flat, the Schr\"odinger-Lichnerowicz formula implies
  that $\ind(D_E)\in K_*(A)$ is an obstruction to positive scalar curvature.
\item On the other hand, the universal property of $C^*_{\max}\pi_1(M)$ implies
  that the representation of $\pi_1(M)$ which gives rise to the bundle $E$
  induces a $C^*$-algebra homomorphism $C^*_{\max}\pi_1(M)\to A$. Moreover,
  the induced map in K-theory sends $\alpha_{\max}(M)\in K_*(C^*_{\max}\pi_1(M))$ to
  $\ind(D_E)\in K_*(A)$.
\item Finally, an index theorem computes $\ind(D_E)$ in terms of the
  degrees of the maps $f_i$ and in particular shows that $\ind(D_E)\ne
  0$. It follows that $\alpha_{\max}(D)\ne 0\in K_*(C^*_{\max}\pi_1(M))$.
\end{enumerate}

A main innovation is the technically quite non-trivial
construction of an 
associated honestly flat bundle, but with infinite dimensional fibers. 

In \cite{HankeKotschickRoeSchick} we relate enlargeability to the classical
strong Novikov conjecture, 
which deals with 
$C^*_{\text{r}}\Gamma$ instead of $C^*_{\max}\Gamma$. 

\begin{figure}
\begin{center}
\epsfig{file=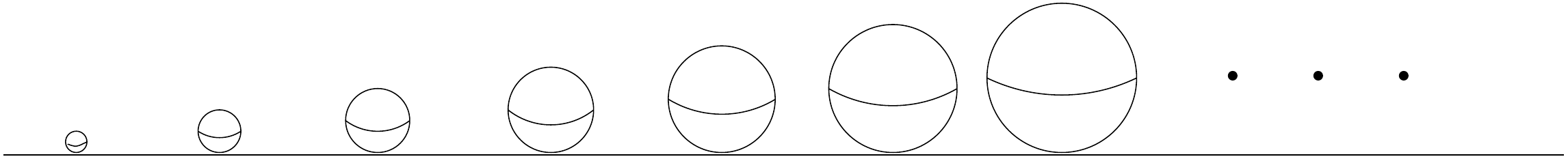,width=13cm}
\caption{The connected balloon space $B_\infty$}\label{fig:connectedballoons}
\end{center}
\end{figure}

The main idea here is to use the functoriality of the large scale index.
The argument becomes technically easier if we assume that all the covering
spaces which determine the enlargeability of $M$ are the universal covering
$\tilde M$. In this situation, the first main point is that the geometry
allows us to combine all the maps $f_i\colon \tilde M\to S^n$ into one map
$F\colon \tilde M\to B_\infty$, where $B_\infty$ is the ``infinite balloon
space'', sketched in Figure~\ref{fig:connectedballoons}.
It is defined using a collection of $n$-spheres of increasing radii $i = 1,2,3, \ldots$, with the sphere of radius $i$ 
attached to the point $i \in [0, \infty)$ at the south pole of $S^n$, and is
equipped with the path metric.  

Using the Mayer-Vietoris sequence and induction on the dimension, one can
calculate the coarse K-homology of $B_\infty$ and the K-theory of its Roe
algebra $C^*B_\infty$. In particular, we obtain $KX_n(B_\infty) \iso
\prod_{i\in\naturals}\integers/\bigoplus_{i\in\naturals}\integers$ and a
direct calculation allows us to establish the coarse Baum-Connes conjecture 
for this space. We obtain a commutative diagram of
K-homology and K-theory groups as follows:
\begin{equation*}
  \begin{CD}
  K_n(M) @>{\iso}>>  K_{n+1}(D^*\tilde M^\Gamma/C^*\tilde M^\Gamma)  @>>>
    K_n(C^*\tilde M^\Gamma) @>{\iso}>> K_n(C^*_{\text{r}}\Gamma)\\
 &&   @VVV @VVV \\
 &&   K_{n+1}(D^*\tilde M/C^*\tilde M) @>>> K_n(C^*\tilde M)\\
 &&  @VVV @VVV\\
 && K_{n+1}(D^*B_\infty/C^*B_\infty)  @>>> K_n(C^*B_\infty)\\
 &&  @VVV   @VV=V\\
 \prod \integers/\bigoplus\integers @>{\iso}>>
 KX_n(B_\infty)   @>{\iso}>> K_n(C^*B_\infty)
   \end{CD}
\end{equation*}
A topological calculation allows one to work out the image
of the fundamental class of $M$ in $KX_n(B_\infty)\iso \prod \integers/\bigoplus
\integers$: it is the class represented by the sequence of degrees
$(\deg(f_i))_{i\in\naturals}$ which by assumption is non-zero. Because of the
coarse Baum-Connes conjecture, the image in the bottom right corner is also
non-zero, which finally implies that also the image $\alpha(M)\in
K_n(C^*_{\text{r}}\Gamma)$ is non-zero, as claimed by the theorem.


\section{Vanishing of the index under partial positivity}

The main reason why one can apply index theory to geometric and topological
questions is that a special geometric situation implies vanishing results for
the index.
It is very important to develop further instances of such vanishing theorems,
in order to widen the scope of the consequences of the index method. Here we 
present one of these, which is valid in the context of large scale index theory:

\begin{theorem}\label{theo:partialvanish}
  Let $M$ be a complete non-compact connected Riemannian spin manifold. Let
  $E\to M$ be a flat bundle of Hilbert $A$-modules for a $C^*$-algebra
  $A$. Assume that the scalar curvature is uniformly positive
  outside a compact subset
  . 

  Then the large scale index of the Dirac operator twisted with $E$ vanishes.
\end{theorem}

For $A=\complexs$, this result has been stated by John Roe 
\cite{Roe_index_coarse, roe2012positive}. A different proof, which covers the
general case, is given by Bernhard Hanke, Daniel Pape and the author in
\cite[Theorem 3.11]{HankePapeSchick}. 

A concrete application of Theorem \ref{theo:partialvanish} to
compact spin manifolds is the codimension-$2$ obstruction to positive scalar
curvature of Theorem \ref{theo:codim2}. In its proof in
\cite{HankePapeSchick}, a gluing and bending 
construction of an intermediate space gives positive scalar curvature outside
of a neighborhood of the hypersurface.



\section{The Stolz exact sequence}
\label{sec:stolz_seq}

Stephan Stolz (compare \cite[Proposition 1.27]{piazzaschickrho}) introduced a
long exact sequence that makes systematic the 
bordism classification of metrics of positive scalar curvature. It is quite
similar in spirit to the surgery exact sequence for the classification of
closed manifolds.

\begin{convention}
  Throughout the remainder of the article, a Riemannian metric on a manifold
  with boundary is 
   assumed to have product structure near the boundary.
\end{convention}

\begin{definition}
  Fix a reference space $X$.
  \begin{enumerate}
  \item The group $\Omega_n^{\spin}(X)$ is the usual spin bordism group,
    consisting of cycles $f\colon M\to X$, with $M$ a closed $n$-dimensional
    spin manifold. The equivalence relation is spin bordism.
  \item The \emph{structure group} $\Pos_n^{\spin}(X)$ is the group of bordism
    classes of metrics of positive scalar curvature on $n$-dimensional closed
    spin manifolds with reference map to $X$. Two such manifolds
    $(M_i,g_i,f_i\colon M_i\to X)$ are called \emph{bordant} if there is a
    manifold $W$ with boundary, with metric $G$ of positive scalar curvature
    and with reference map $F\colon W\to X$ such that its boundary is
    $M_1\disjointunion (-M_2)$ and $G,F$ restrict to the given $g_i,f_i$ at the
    boundary. 
  \item Finally, the group $R_n(X)$ is the group of equivalence classes of
    compact spin manifolds $W$ 
    with boundary, with reference map $f\colon W\to X$ and with a metric $g$
    of positive scalar curvature on $\boundary W$. Again, the equivalence
    relation on such cycles is bordism, where a bordism between
    $(W_1,f_1,g_1)$ and $(W_2,f_2,g_2)$ is a manifold $Y$ with boundary
    $\boundary Y= W_1\cup_{\boundary W_1} Z \cup_{\boundary W_2} -W_2$ (where
    $Z$ is a spin bordism 
    between $W_1$ and $W_2$) together with a continuous map $f\colon Y\to X$
    and a positive scalar curvature metric $g$ on $Z$. Of course, the
    restriction of $f$ to $W_j$ must be $f_j$ and the restriction of $g$ to
    $\boundary W_j$ must be $g_j$. It turns out that $R_n(X)$ only depends on $\pi_1(X)$.
  \item The group structure in each of the three cases is given by disjoint
    union, and the inverse is obtained by reversing the spin structure and
    leaving all other data unchanged.
  \item There are evident ``forget structure'' and ``take boundary'' maps
    between these groups. Using these, one obtains a long exact sequence, the
    \emph{Stolz positive scalar curvature exact sequence}
    \begin{equation}
      \label{eq:Stolz}
      \cdots  \to R_{n+1}(\pi_1(X)) \to \Pos^{\spin}_n(X) \to
      \Omega^{\spin}_n(X)\to \cdots 
    \end{equation}
  \end{enumerate}
\end{definition}

The most useful cases of this sequence arises if $X=M$ and $f=\id\colon M\to
M$ or if $X=B\Gamma$ is the classifying
space of a discrete group and $f\colon M\to B\Gamma$ induces the identity on
the fundamental groups.

To get information about $\Pos^{\spin}_n(M)$ we want to use index theory
systematically by mapping in a consistent way to the analytic exact sequence
of Nigel Higson and 
John Roe. This sequence is simply the long exact K-theory sequence of the
extension
$  0\to C^*\tilde M^\Gamma \to D^*\tilde M^\Gamma\to
  D^*\tilde M^\Gamma/C^*\tilde M^\Gamma \to 0$, where $\tilde M$ is the
  universal covering of $M$, and $\Gamma=\pi_1(M)$.

Using that $K_*(C^*\tilde M^\Gamma)=K_*(C^*_{\text{r}}\Gamma)$ and that
$K_{*+1}(D^*\tilde M^\Gamma/C^*\tilde M^\Gamma) = K_*(M)$, we obtain the following
theorem, compare \cite[Theorem 1.39]{piazzaschickrho}
\begin{theorem} Let $X$ be a topological space with $\Gamma$-covering $\tilde
  X$.
We have a natural canonical commutative diagram (established with full proof
if $n$ is odd)
\begin{equation*}
{\small  \begin{CD}
   @>>> \Omega^{\spin}_{n+1} (X) @>>>   R_{n+1}(X) @>>>
    \Pos^{\spin}_n (X) @>>> \Omega^{\spin}_n (X) \cdots\\ 
    && @VV{\beta}V @VV{\ind}V  @VV{\rho_\Gamma}V @VV{\beta}V\\
     @>>>  K_{n+1} (X) @>>> K_{n+1} ( C^*_{\text{r}}\Gamma) @>>>
     K_{n+1}(D^*\tilde X^\Gamma)  @>>>  K_{n} ( X) \cdots\\
    \end{CD}
}
\end{equation*}

Here, $\beta$ is obtained by taking the large scale (equivariant) index of the
Dirac operator on the covering of a cycle $f\colon N\to X$ and then
using 
functoriality of large scale index theory to push forward via $f_*$ from $K_*(C^*\tilde
N^\Gamma)$ to $K_*(C^*\tilde X^\Gamma)$. It coincides with the Atiyah
orientation as natural transformation from spin bordism to K-homology.
 Similarly, $\rho_\Gamma$ is obtained by
constructing the structure invariant of the positive scalar curvature metric
of the covering of $(N,g,f\colon N\to X)$ and then using naturality to
push forward along $f_*$ from $K_*(D^*\tilde N^\Gamma)$ to $K_*(D^*\tilde
X^\Gamma)$. 

Finally, $\ind$ assigns to a manifold with boundary and positive scalar
curvature at the boundary an Atiyah-Patodi-Singer type index.
\end{theorem}

Note that the assertion of Theorem \ref{theo:rhoAPS} is that the (index based)
maps all exist, that they are indeed well defined, i.e.~invariant under
bordism, and that the diagram is commutative.
This means that we have to work (for the cycles and for the equivalence
relation) throughout with manifolds with boundary. It turns out that large
scale index theory can very elegantly and efficiently be used to carry out
index theory on manifolds with boundary, as well.

\section{Index theory on manifolds with boundary}
\label{sec:index_boundary}

Our method to do index theory on a manifold with boundary is 
simply to  attach an infinite half-cylinder to the boundary. This produces a
manifold without boundary, of course at the expense that the resulting
manifold is never compact. However, large scale index theory can deal with such
spaces. 

To obtain the appropriate information, the construction must take the extra
information into account coming from the fact that the metric of the boundary
is assumed to have positive scalar curvature. Let us review the construction:
\begin{enumerate}
\item We start with a smooth manifold $W$ with boundary, with a Riemannian
  metric $g$ which has positive scalar curvature near the boundary (and a
  product structure there, by our general convention). As a metric space, $W$
  is assumed to be complete. Moreover, we fix an
  auxiliary $C^*$-algebra $A$ and a flat
  Hilbert $A$-module bundle $E$ on $W$ (again with product structure near the
  boundary).
\item We now attach a half-cylinder $\boundary W\times [0,\infty)$ to the
  boundary to obtain $W_\infty$ and extend all the structures over
  $W_\infty$. We obtain a complete manifold without boundary, with product end
  $\boundary W\times [0,\infty)$.
\item As in Section \ref{sec:subs-theory-manif}, the Dirac operator $D_E$
  produces a bounded operator $\psi(D_E)$ in $D^*(W_\infty;A)$.
\item Now, however, we use the invertibility of $D_E$ on $\boundary M\times
  [0,\infty)$ coming from the Schr\"odinger-Lichnerowicz formula:
for suitable $\psi$ the element  $1-\psi(D_E)^2$ does not only lie in
$C^*(W_\infty;A)$ but in the smaller
  ideal $C^*(W\subset W_\infty;A)$. This is by definition generated by all
  locally compact finite propagation operators $T$ which are \emph{supported
    near $W$}, which means that there is $R>0$ such that $T\phi=0$ and $\phi
  T=0$ 
  whenever $\phi$ is a compactly supported function with $d(\supp(\phi),W)>R$.
\item Correspondingly, the fundamental class of $W_\infty$ has a canonical
  lift to a class $[W,g|_{\boundary W}]$ in
  $K_{n+1}(D^*(W_\infty;A)/C^*(W\subset W_\infty;A))$. This class does in
  general depend on the positive scalar curvature metric on the boundary.
\item As usual, we next define the ``large scale Atiyah-Patodi-Singer index''
  by applying the
  boundary map of the long exact K-theory sequence, now for the extension
  \begin{equation*}
    0\to C^*(W\subset W_\infty;A) \to D^*(W_\infty;A)\to
    D^*(W_\infty;A)/C^*(W\subset W_\infty;A)\to 0,
  \end{equation*}
  to obtain $\ind(D_W, g_{\boundary W})\in K_n(C^*(W\subset W_\infty;A))\iso
  K_n(C^*(W;A))$. The 
  latter isomorphism is induced by the inclusion $ C^*(W;A)\into
  C^*(W\subset 
  W_\infty;A)$ which just extends the operators by zero. Note that this
  construction of the index required the invertibility of the operator at the
  boundary and indeed depends in general on the metric of positive scalar
  curvature at $\boundary W$.
\end{enumerate}

Note that large scale index theory in the situation we just described produces
two invariants which depend on the positive scalar curvature metric on
$\boundary W$, namely $\ind(D_W,g_{\boundary W})\in K_n(C^*(W\subset
W_\infty;A))$, but also the secondary invariant $\rho(\boundary W, g_{\boundary
  W}) \in K_n(D^*(\boundary W;A))$. A major result, which we consider a
secondary higher Atiyah-Patodi-Singer index theorem, relates these two.

\begin{theorem}\label{theo:rhoAPS} (compare \cite[Theorem
  1.22]{piazzaschickrho})
  Let $(W,g_W)$ be an even dimensional Riemannian spin-manifold with boundary
$\boundary W$
such that $g_{\boundary W}$ 
has positive scalar curvature. Assume that a group $\Gamma$ acts
isometrically on $M$. Then 
\begin{equation*}
\iota_* (\ind_\Gamma (D_W))= j_*(\rho(\boundary W, g_{\boundary W})) \quad\text{in}\quad K_{0}
(D^*(W)^\Gamma).
\end{equation*}
Here, we use $j\colon D^*(\boundary W)^\Gamma\to D^*W^\Gamma$ induced by the
inclusion 
$\boundary W\to W$ and $\iota\colon C^*(W)^\Gamma\to D^*(W)^\Gamma$ the inclusion.
\end{theorem}

\begin{remark}
  Above we apply the obvious generalization of the construction of the large
  scale index of Sections \ref{sec:subs-theory-manif} and
  \ref{sec:index_boundary} to an equivariant situation, where subalgebras
  $D^*(W)^\Gamma$ and $C^*(W)^\Gamma$ generated by invariant operators are
  used.
  This works because the Dirac operator then itself is invariant under the
  group $\Gamma$.
\end{remark}

\begin{remark}
  The heart of the proof of Theorem \ref{theo:rhoAPS} is an explicit secondary
  index calculation in a model (product) case which is surprisingly
  intricate. 
\end{remark}

\begin{remark}
  The assertion of Theorem \ref{theo:rhoAPS} should generalize to arbitrary
  (non-cocompact) spin manifolds, to Hilbert $C^*$-algebra coefficient bundles
  and
  to the K-theory of real $C^*$-algebras. In a preprint of Zhizhang Xie and
  Guoliang Yu \cite{XieYu} an argument is sketched which shows how to extend
  the result to arbitrary dimensions and to real $C^*$-algebras. 
\end{remark}

\section{Constructions of new classes of metrics of positive scalar curvature}
\label{sec:HankeSteimle}

The fundamental idea in the construction of the ``geometrically significant''
homotopy classes of the space of metrics of positive scalar curvature of
Theorem \ref{theo:large_pi} is quite old and based on index
theory:

Given a closed $n$-dimensional spin manifold $B$ with $\hat A(B)\ne
0$, we know that $B$ does not admit a metric of positive scalar curvature.

Remove an embedded disc from $B$. The result is a manifold $W$
with boundary $\boundary W=S^{n-1}$. Given any metric of positive scalar
curvature 
on $W$ (with product structure near the boundary), the corresponding
boundary metric $g_1$ can not be homotopic to the standard metric on $S^{n-1}$
because then 
one could glue in the standard 
disc (with positive scalar curvature) to obtain a metric of positive scalar
curvature on $B$. Now, if $g_1$ is homotopic to $\psi^*g_{\eucl}$ for a
non-identity diffeomorphism we can glue the disc back in with $\psi$ to obtain
a new manifold $B_\psi$ which is of positive scalar curvature. Note that
$B_\psi$ is not necessarily diffeomorphic to $B$, but using the Alexander
trick there is a homeomorphism between $B_\psi$ and $B$. As the rational
Pontryagin classes and therefore the $\hat A$-genus are homeomorphism
invariant, $\hat A(B_\psi)\ne 0$, also $B_\psi$ can not carry a
metric of positive scalar curvature.

Observe that exactly the same argument can be applied to a family situation:
Let $Y\to S^k$ be a family (i.e.~bundle) of manifolds with boundary, with
boundary $S^n\times S^k$. Assume there is a family of metrics $g_x$ of positive
scalar curvature 
on $Y$ (product near the boundary). If this family of metrics is homotopic to
the constant family 
consisting of the standard metric (or to a pullback of that one along a family
of diffeomorphisms $\psi_x$) we can glue in $D^{n+1}\times S^k$ to obtain a
family of closed manifolds which admits a metric of positive scalar curvature
(fiberwise, and then also the total space $X$ admits such a metric).
Note that this is only interesting if each $g_x$ is in the component of the
standard metric, which we therefore assume.

Alas: if the total space $X$ has non-trivial $\hat A$-genus, this is a
contradiction (and again, by the homeomorphism invariance and the Alexander
trick the argument works modulo diffeomorphism).

Note that this requires two important ingredients:
\begin{enumerate}
\item the topological situation with the bundle $Y$ (and $X$)
\item  the geometric input of a family of metrics of
  positive scalar curvature on $Y$.
\end{enumerate}

It turns out that already the topological input is surprisingly
difficult to get. It means that (after the gluing) we have a fiber bundle
$M\to X\to 
S^k$ of spin manifolds where $M$ does admit a metric of positive scalar
curvature, therefore $\hat A(M)=0$, whereas $\hat A(X)\ne 0$. Note that this
means that the $\hat A$-genus is not multiplicative in fiber bundles, even if
the base is simply connected (in contrast to the L-genus). 

In \cite[Theorem 1.4]{HSS} we use advanced differential topology, in
particular surgery theory, Casson's theory of prefibrations and Hatcher's
theory of concordance spaces to prove that the required fiber bundles $X$ exist:
\begin{theorem}\label{theo:Ahatbundles}
  For sufficiently large $n$, there are $4n$-dimensional smooth closed spin
  manifolds $X$ with non-vanishing $\hat A$-genus fitting into a smooth fiber
  bundle $F\to X\to S^k$ such that $F$ admits a metric of positive scalar
  curvature, is highly connected and the bundle contains as subbundle
  $D^{4n-k}\times S^k$.
\end{theorem}

How about the second ingredient, the existence of the family of metrics of
positive scalar curvature on $Y=X\setminus D^{4n-k}\times S^k$?

The only tool known which can provide such metrics is the surgery method of
Gromov and Lawson. In highly non-trivial work \cite{Walsh2} this has
been extended by Mark Walsh to families of the kind $X$ as constructed in
Theorem \ref{theo:Ahatbundles}. To apply this, we use the high connectivity
and results of Kiyoshi Igusa on Morse theory for fiber bundles
\cite{Igusa}. As a consequence we obtain 
Theorem \ref{theo:large_pi}.

\section{Open problems}

The geometry of positive scalar curvature and the development and application
of large scale index theory is a vibrant field of research, with a host of
important open problems. Many of those were already mentioned above; here we
want to highlight them and add a couple of further
directions of research.

\bigskip\noindent\textbf{Gromov-Lawson-Rosenberg conjecture}.
We should find further obstructions to positive scalar curvature on spin
manifolds, in particular for
finite fundamental group $(\integers/p\integers)^n$ for the so-called toral manifolds. We expect that this will
require fundamentally new ideas.

On the other hand, can the class of fundamental groups for which the
conjecture holds be 
described systematically?

\bigskip\noindent
\textbf{Stable Gromov-Lawson-Rosenberg conjecture}.
The stable Gromov-Lawson-Rosenberg conjecture and its weaker cousin
\ref{conj:Rosenberg_power} which states that ``the Rosenberg index sees
everything about positive scalar curvature which can be seen using the Dirac
operator'' follow from the strong Novikov conjecture. It would be spectacular
to find counterexamples to either of these (they are expected to exist, but to
find them will of course be very hard). It is necessary to investigate this
question further. In this context, the role of the codimension-$2$ obstruction
as discussed in Theorem \ref{theo:codim2} should be understood. 

This theorem  should extend to the signature operator, which
will shed new light on its meaning. Vaguely, we conjecture
the following.

\begin{conjecture}\label{conj:codim2sig}
  Let $M_1,M_2$ be two complete non-compact connected oriented Riemannian
  manifolds and $f\colon M_1\to M_2$ a sufficiently well behaved map which is
  an ``oriented homotopy equivalence near infinity''.

  Let
  $E\to M_2$ be a flat bundle of Hilbert $A$-modules for a $C^*$-algebra
  $A$. Let $D^{\sgn}_E$ be the signature operator on $M_1$ twisted with the
  flat bundle $E$,
  and $D^{\sgn}_{f^*E}$ the signature operator on $M_2$ twisted by $f^*E$. Then
  the large scale indices of these two operators should coincide,
  i.e.
  $$f_*(\ind(D^{\sgn}_{f^*E})) = \ind(D^{\sgn}_E) \in K_*(C^*(M_1;A)).$$
\end{conjecture}

Note that, in this conjecture, one has to work out the
  precise concept of ``sufficiently well behaved'' and of ``homotopy
  equivalence at infinity''.

\bigskip
\noindent
\textbf{Area based large scale geometry}.
Large scale geometry is based on the metric spaces and distances, viewed from
a coarse perspective. Curvature, on the other hand, is a concept based on the
bending of surfaces, where scalar curvature looks at the average over all
possible surface curvatures through a given point.

This is reflected in the fact that area-enlargeability suffices to obstruct
positive scalar curvature (Theorem \ref{theo:enlarg}). So far, this
is not captured well by large scale index theory. 

This suggests that a program should be developed for large scale geometry
based on $2$-dimensional areas instead of lengths. A possible starting point
would be to work on a relative of the loop space and carry out the analysis
there. This is interesting in its own right, with a host of potential further
applications, but seems to require new analytical tools.

Note that the axiomatic abstraction from metric spaces to coarse structures as
developed by John Roe \cite{RoeCoarseGeom} does not seem to apply here. Of
course, this generalization is interesting in its own
right and applications to positive scalar curvature should be developed
further. 

\bigskip\noindent
\textbf{Coarse Baum-Connes conjecture}.
If the coarse Baum-Connes conjecture holds for the classifying space of a
discrete group, then also the strong Novikov conjecture is true for this
group. Moreover, the validity of the coarse Baum-Connes conjecture is a
powerful computational tool. On the other hand, there are enigmatic
counterexamples due to Guoliang Yu if one drops the ``bounded geometry''
condition on the space in question. 

We expect that many new classes of metric spaces can be found where the coarse
Baum-Connes conjecture can be established. But we also feel that the search
for counterexamples (of bounded geometry) should be intensified.

\bigskip\noindent
\textbf{Aspherical manifolds}.
A lot of attention has been given to the special class of aspherical manifolds.

\begin{conjecture}\label{conj:aspherical}
  Let $M$ be a closed manifold whose universal covering is contractible
  (i.e.~$M$ is \emph{aspherical}). Then $M$ does not admit a metric of
  positive scalar curvature.
\end{conjecture}
Often the geometry implies that a manifold is aspherical
(e.g.~if it admits a metric which is non-positively curved in the sense of
comparison geometry). The conjecture states that in a weak sense this is the
only way a manifold can be aspherical. The strong Novikov conjecture for an
aspherical spin manifold $M$ implies that $\alpha(M)\ne
0$ because the Dirac operator of a manifold $M$ always represents a non-zero
K-homology class in $K_*(M)$, and here $M=B\pi_1(M)$. Of course, we now look
for ways to directly use the asphericity in proofs of (special cases) of
Conjecture \ref{conj:aspherical}.

\bigskip\noindent
\textbf{Mapping surgery to analysis to homology}.
The program to map surgery to analysis has been fully carried out in
\cite{piazzaschickrho} only for half the dimensions, and only for complex
$C^*$-algebras, based on a delicate explicit index calculation. New
developments, in particular the work of Zhizhang Xie and Guoliang Yu
\cite{XieYu} extend this to all dimensions with a modified method. It remains
to develop the details of this (or an alternative) approach and to carry it
over to the more powerful real $C^*$-algebras.

K-theory of $C^*$-algebras is a very powerful tool. Most useful, however is
its combination with homological tools (in particular Hochschild and cyclic 
(co)homology). To achieve this systematically and use it for the
classification of metrics of positive scalar curvature, we propose a program
to not only map the positive scalar curvature sequence to analysis, as
described in Section \ref{sec:stolz_seq}, but then to map further to a
corresponding long exact sequence of (cyclic) homology groups. There, one
would then see primary and secondary numerical invariants of higher index
theory.

The primary invariants are well developed. Rather not understood,
however, is the theory related to the 
secondary invariants. Indeed, the relevant algebra $D^*M$ is very large and
the usual dense and holomorphically closed subalgebras on which explicit
formulas for the Connes-Karoubi Chern character would make sense seem hard to
come by. 
Exactly because of this we feel that the development of such a theory will
shed important new light on the power of the secondary invariants used in
Section \ref{sec:stolz_seq}.

Apart from this general program, the theory described above needs to be
applied in more concrete situations. The analytic structure set
 $K_*(D^*M)$, despite its evident potential, so far has only
been used in a small number of concrete contexts (compare in particular
\cite{HR4}). This must change before we will have a definite idea of its
power.

\bigskip\noindent\textbf{Minimal surface obstructions to positive scalar
  curvature}.
The minimal surface method is the only known tool to obstruct the existence of
a metric of positive scalar curvature which works for non-spin manifolds of
dimension $\ge 5$. So far, it is controversial how to extend the method
if the dimension of the manifold in question is larger than $8$ (due to
singularities which develop in the minimal hypersurfaces one wants to
use). Joachim Lohkamp has a program to achieve this.

The minimal surface technique has so far only been used together with
the Gauss-Bonnet theorem (via an iterative approach, until the hypersurface
is $2$-dimensional). Are there other ways to exploit it and combine it with
obstructions to positive scalar curvature (Seiberg-Witten, spin Dirac) and
what are the relations?

\bigskip\noindent\textbf{Enlargeability and non-spin manifolds}.
As just one case of the question whether the known results for spin manifolds
carry over to non-spin manifolds consider the following:
\begin{question}
  Let $M$ be an arbitrary (non-spin) closed $n$-dimensional manifold with
  coverings
  $M_\epsilon$ and with  maps $f_\epsilon\colon M_\epsilon
  \to S^n$ which are constant outside a compact subset of $M_\epsilon$, which
  have non-zero degree and which are $\epsilon$-contracting.

  Can $M$ admit a metric of positive scalar curvature?
\end{question}
Using the minimal surface method, Gromov and Lawson \cite[Section
12]{GromovLawson} have shown that this is not the case if $n\le 7$.



\bigskip\noindent\textbf{Topology of the space of metrics of positive scalar curvature}.
The stable Gromov-Lawson-Rosenberg conjecture shows that there is a stability
feature in the topology of the space $\Riemplus(M)$ of metrics of positive
scalar curvature: if index theory suggests that it should be non-empty, this
might be violated by $M$ itself, but after iterated product with $B$,
eventually it is 
non-empty. Are there similar stability features concerning the (higher)
homotopy groups of $\Riemplus(M)$? It would also be important to develop
estimates on the stable range.

\bigskip\noindent\textbf{The space of metrics of positive scalar curvature and fundamental group}.
We know well that, for a spin manifold $M$ with a complicated fundamental group
$\Gamma$, the existence of a metric of positive scalar curvature is not only
obstructed by $\hat A(M)$, but by
 $\alpha(M)\in K_*(C^*\Gamma)$, and many elements of $K_*(C^*\Gamma)$ are
 indeed realized as values of $\alpha(M)$.

Similarly, we should expect that the topology of $\Riemplus(M)$, if non-empty,
should be governed by $K_*(C^*\Gamma)$. At the moment, a precise conjecture
(e.g.~about the homotopy groups) seems too far-fetched. Still, the methods of
large scale and higher index theory should be developed to the point that they
are available to detect new elements in $\pi_*(\Riemplus(M))$, and
one should systematically construct non-trivial examples.




\bibliographystyle{plain}
\bibliography{icm}

\end{document}